\documentclass[11pt]{article}
\usepackage{amsfonts}
\usepackage{amsmath}

\input{tcilatex}
\begin{document}

\title{Almost absolute weighted summability with index k}
\author{Mehmet Ali Sar\i g\"{o}l$^{1}$ and Mohammad Mursaleen$^{2}$ \\
1 {\small {Department of Mathematics, Faculty of Art and Science, Pamukkale
University, Denizli 20160, Turkey}}\\
msarigol{\small {@pau.du.tr} }\\
2 Department of Mathematics, Aligarh Muslim University, Aligarh 202 002, Up
India\\
mursaleenm@gmail.com}
\maketitle

\begin{abstract}
The space $\widehat{\ell }_{k}$ of absolutely almost convergent series was
introduced and studied by Das et al [$4$]$,$ which plays an important role
in summability theory, approximation theory, Fourier analysis, etc. In the
present paper we generalize the space making use of some factors and
weighted mean transformations, investigate its toplogical structures and
relations between classical sequence spaces. Also we characterize certain
matrix operatos on it.
\end{abstract}

\section{\textbf{Introduction}}

Any subspace of $w,$ the set of all sequences of complex numbers, is called
a \textit{sequence space}. A $BK$ space $X$ is a Banach sequence space with
the property that the map $p_{n}:X\rightarrow 
\mathbb{C}
$ defined by $p_{n}(x)=x_{n}$ is continuous for all $n\geq 0,$ where $%
\mathbb{C}
$ denotes the complex field. Let $\ell _{\infty }$ be the subspace of all
bounded sequences of $w.$ A sequence $\left( x_{n}\right) \in \ell _{\infty
} $ is said to be almost convergent to $\gamma $ if all of its Banach limits 
$\left[ 1\right] $ are equal to $\gamma $. Lorentz $\left[ 8\right] $
characterized almost convergence that a sequence $\left( x_{n}\right) $ is
almost convergent to $\gamma $ if and

{\small \textit{2010 AMS Subject Classification: }}40C05, 40D25, 40F05, 46A45

{\small \textit{Keywords: Sequence space, Bk-space, Absolute and almost
summability, Matrix transformations}}\newline
only if 
\begin{equation}
\frac{1}{m+1}\tsum_{v=0}^{m}x_{n+v}\rightarrow \gamma \text{ as }%
m\rightarrow \infty \text{ uniformly in }n  \tag{1.1}
\end{equation}

This notation plays an important role in summability theory, approximation
theory and Fourier analysis and was investigated by several authors. For
example, it was later used to define and study some concepts such as
conservative and regular matrices, some sequence spaces and matrix
transformations (see [2], [6], [7], [9], [10], [11], [15]).

Absolute almost convergence emerges naturally as absolute analogue of almost
convergence just as absolute convergence emerged out of the concept of
convergence. To introduce this concept, let $\Sigma a_{v}$ be a given
infinite series with $s_{n}$ as its n-th partial sum. The series $\Sigma
a_{v}$ is said to be absolutely almost convergent series if (see [3]) 
\begin{equation*}
\tsum_{m=0}^{\infty }\left\vert \psi _{m,n}\right\vert ^{k}<\infty ,\text{ }%
k>0,
\end{equation*}%
uniformly in $n,$ where 
\begin{equation*}
\psi _{m,n}=\left\{ 
\begin{array}{c}
a_{n},\text{ \ }m=0 \\ 
\frac{1}{m(m+1)}\tsum_{v=1}^{m}va_{n+v},\text{ }m\geq 1.%
\end{array}%
\right.
\end{equation*}%
The space of all absolutely almost convergent series 
\begin{equation*}
\widehat{\ell }_{k}=\left\{ a:\tsum_{m=0}^{\infty }\left\vert \psi
_{m,n}\right\vert ^{k}<\infty ,\text{ uniformly in }n,\text{ }k>0\right\} .
\end{equation*}%
was first defined and studied in $\left[ 4\right] $. We note an important
relation between $\widehat{\ell }$ and absolute Cesaro summability $%
\left\vert C,1\right\vert $ in Flett's notation [5], $\widehat{\ell }\subset 
$ $\left\vert C,1\right\vert $ $\left[ 4\right] .$

\section{Main Results}

The purpose of the present paper is to define an absolute almost weighted
summability using some factors and weighted means, which extends the well
known concept of absolute almost convergence of Das et al $\left[ 4\right] ,$
and to study its topological structures. Also we investigate relations
between classical sequence spaces and characterize certain matrix operatos
on it.

For any sequence $\left( s_{n}\right) ,$ we define $T_{m,n}$ by 
\begin{equation*}
T_{-1,n}(s)=s_{n-1},\text{ }T_{m,n}(s)=\frac{1}{P_{m}}%
\tsum_{v=0}^{m}p_{v}s_{n+v},m\geq 0.
\end{equation*}%
A straightforward calculation then shows that%
\begin{equation*}
F_{m,n}(a)=T_{m,n}(s)-T_{m-1,n}(s)=\left\{ 
\begin{array}{c}
a_{n},\text{ \ }m=0 \\ 
\frac{p_{m}}{P_{m}P_{m-1}}\tsum_{v=1}^{m}P_{v-1}a_{n+v},\text{ }m\geq 1,%
\end{array}%
\right.
\end{equation*}%
where $(p_{n})$ is a sequence of positive real numbers with 
\begin{equation}
P_{n}=p_{0}+p_{1}+...+p_{n}\rightarrow \infty \ as\ n\rightarrow \infty ,\
P_{-1}=p_{-1}=0.  \tag{2.1}
\end{equation}%
So we give the following definition.

\textbf{Definition 2.1. }Let $\Sigma a_{v}$ be an infinite series with
partial summations $s_{n}.$ Let $(p_{n})$ and $\left( u_{n}\right) $ be
sequences of positive real numbers. The series $\Sigma a_{v}$ is said to be
absolute almost weighted summable $\left\vert f(\overline{N}%
_{p}),u_{m}\right\vert _{k},$ $k\geq 1,$\ if 
\begin{equation}
\sum_{m=0}^{\infty }u_{m}^{k-1}\left\vert F_{m,n}(a)\right\vert ^{k}<\infty 
\tag{2.2}
\end{equation}%
uniformly in $n.$

For $\left\vert f(\overline{N}_{p}^{u})\right\vert _{k},k\geq 1,$ we write
the set of all series summable by the method $\left\vert f(\overline{N}%
_{p}),u_{m}\right\vert _{k}.$ Then, $\Sigma a_{v}$ is summable $\left\vert f(%
\overline{N}_{p}),u_{m}\right\vert _{k}$\ iff the series $\Sigma a_{v}\in $ $%
\left\vert f(\overline{N}_{p}^{u})\right\vert _{k}.$ Note that, in the case$%
\ u_{m}=p_{m}=1$ for $m\geq 0,$ it reduces to the set of absolutely almost
convergent series $\widehat{\ell }_{k}$ given by Das, Kuttner and Nanda $%
\left[ 3\right] .$ Further, it is clear that the space $\left\vert \overline{%
N}_{p}^{u}\right\vert _{k}$ is derived from $\left\vert f(\overline{N}%
_{p}^{u})\right\vert _{k}$ by putting $n=0$ (see [9], [12], [13]), and also $%
\left\vert f(\overline{N}_{p}^{u})\right\vert _{k}\subset \left\vert 
\overline{N}_{p}^{u}\right\vert _{k},$ but the converse is not true.

First it would be appropriate to clarify some relations between the new
method and classical sequence spaces such as $bs$ and $\ell _{\infty },$
where $bs$ and $\ell _{\infty }$ are the set of all bounded series and
sequences, respectively.\smallskip

\textbf{Theorem 2.2. }Let $(p_{m})$ and $(u_{m})$\ be two sequences of
positive numbers. If

$(i)$ If 
\begin{equation}
\left( \frac{1}{u_{m}}\right) \in \ell _{\infty },  \tag{2.3}
\end{equation}%
then $\left\vert f(\overline{N}_{p}^{u})\right\vert _{k}\subset \ell
_{\infty },$ $k\geq 1.$

$\left( ii\right) $ If 
\begin{equation}
\sum_{m=0}^{\infty }u_{m}^{k-1}\left( \frac{p_{m}}{P_{m}}\right) ^{k}<\infty
,  \tag{2.4}
\end{equation}%
then $bs\subset \left\vert f(\overline{N}_{p}^{u})\right\vert _{k},$ $k>1.$

\textbf{Proof. }$\left( i\right) $ Given $a=(a_{v})\in $ $\left\vert f(%
\overline{N}_{p}^{u})\right\vert _{k}.$ Then, by the definition, there
exists an integer $M$ such that, for all $n,$ 
\begin{equation}
\sum_{m=M}^{\infty }u_{m}^{k-1}\left\vert F_{m,n}(a)\right\vert ^{k}\leq 1. 
\tag{2.5}
\end{equation}%
So, by $\left( 2.3\right) $ and $\left( 2.5\right) ,$ we have $\left\vert
F_{m,n}(a)\right\vert \leq u_{m}^{-1+1/k}$ for $m\geq M$ and all $n.$ On the
other hand, for $m\geq 1,$ 
\begin{equation}
a_{m+n}=\frac{P_{m}}{p_{m}}F_{m,n}(a)-\frac{P_{m-2}}{p_{m-1}}F_{m-1,n}(a), 
\tag{2.6}
\end{equation}%
which gives that $a=(a_{v})\in \ell _{\infty }.$ This completes the proof.

$\left( ii\right) $ Let $a=(a_{v})\in bs.$ Say $M=\sup_{v}\left\vert
\tsum_{j=0}^{v}a_{j}\right\vert .$ Then we have 
\begin{eqnarray*}
\left\vert \tsum_{v=1}^{m}P_{v-1}a_{n+v}\right\vert &=&\left\vert
\tsum_{v=1}^{m-1}(-p_{v})\tsum_{j=1}^{v}a_{n+j}+P_{m-1}%
\tsum_{j=1}^{m}a_{n+j}\right\vert \\
&\leq &\tsum_{v=1}^{m-1}p_{v}\left\vert \tsum_{j=n+1}^{n+v}a_{j}\right\vert
+P_{m-1}\left\vert \tsum_{j=n+1}^{n+m}a_{n+j}\right\vert \\
&\leq &\tsum_{v=1}^{m-1}p_{v}\left\vert
\tsum_{j=0}^{n+v}a_{j}-\tsum_{j=0}^{n}a_{j}\right\vert +P_{m-1}\left\vert
\tsum_{j=0}^{n+m}a_{n+j}-\tsum_{j=0}^{n}a_{j}\right\vert \\
&\leq &4MP_{m-1}
\end{eqnarray*}%
which implies that

\begin{eqnarray*}
\sum_{m=0}^{\infty }u_{m}^{k-1}\left\vert F_{m,n}(a)\right\vert ^{k}
&=&u_{0}^{k-1}\left\vert a_{n}\right\vert ^{k}+\sum_{m=1}^{\infty
}u_{m}^{k-1}\left\vert \frac{p_{m}}{P_{m}P_{m-1}}%
\tsum_{v=1}^{m}P_{v-1}a_{n+v}\right\vert ^{k} \\
&\leq &(4M)^{k}\sum_{m=0}^{\infty }u_{m}^{k-1}\left( \frac{p_{m}}{P_{m}}%
\right) ^{k}<\infty .
\end{eqnarray*}%
This completes the proof.

For the special case $u_{m}=p_{m}=1$ for all $m\geq 0,$ $\left\vert f(%
\overline{N}_{p}^{u})\right\vert _{k}=\widehat{\ell }_{k}$ and $\left(
2.4\right) $ reduces to $\sum_{m=0}^{\infty }\left( m+1\right) ^{-k}<\infty
. $ So we have the following result in $\left[ 3\right] .$

\textbf{Corollary 2.3. }For $k>1,$ $bs\subset \widehat{\ell }_{k}.$\smallskip

\textbf{Theorem 2.4. }Let $(u_{m})$\ be any sequence of positive numbers.
Then, $\left\vert f(\overline{N}_{p}^{u})\right\vert _{k},$ $k\geq 1,$ is a $%
BK$ space with respect to the norm 
\begin{equation}
\left\Vert a\right\Vert _{\left\vert f(\overline{N}_{p}^{u})\right\vert
_{k}}=\sup \left\{ \sum_{m=0}^{\infty }u_{m}^{k-1}\left\vert
F_{m,n}(a)\right\vert ^{k}:n\in 
\mathbb{N}
\right\} ^{1/k}.  \tag{2.7}
\end{equation}

\textbf{Proof. } It is routine to show that $\left( 2.7\right) $ satisfies
the norm conditions. We only note that $\left( 2.7\right) $ is well defined.
In fact, if $a\in \left\vert f(\overline{N}_{p}^{u})\right\vert _{k}$, then,
as in the proof of part $\left( ii\right) $ of Theorem 2.2, there exists an
integer $M$ such that, for all $n,$ 
\begin{equation*}
\sum_{m=M}^{\infty }u_{m}^{k-1}\left\vert F_{m,n}(a)\right\vert ^{k}\leq 1.
\end{equation*}%
and \ $\left( u_{m}^{1-1/k}\left\vert F_{m,n}(a)\right\vert \right) $ is
bounded for all $m,n\geq 0.$ This gives 
\begin{eqnarray*}
&&\sup \left\{ \sum_{m=0}^{\infty }u_{m}^{k-1}\left\vert
F_{m,n}(a)\right\vert ^{k}:n\in 
\mathbb{N}
\right\} \\
&\leq &1+\sup \left\{ \sum_{m=0}^{M}u_{m}^{k-1}\left\vert
F_{m,n}(a)\right\vert ^{k}:n,m\in 
\mathbb{N}
\right\} <\infty .
\end{eqnarray*}%
To prove that it is a Banach space, let us take arbitrary Cauchy sequence $%
\left( a^{m}\right) ,$ where $a^{m}=\left( a_{v}^{m}\right) \in \left\vert f(%
\overline{N}_{p}^{u})\right\vert _{k}$ for $m\geq 0.$ Given $\varepsilon >0.$
Then there exists an integer $m_{0}$ such that $\left\Vert
a^{m_{1}}-a^{m_{2}}\right\Vert _{\left\vert f(\overline{N}%
_{p}^{u})\right\vert _{k}}<\varepsilon $ for $m_{1},m_{2}>m_{0},$ or,
equivalently, 
\begin{equation}
\sup_{n}\left\{ u_{0}^{k-1}\left\vert a_{n}^{m_{1}}-a_{n}^{m_{2}}\right\vert
^{k}+\sum_{m=1}^{\infty }u_{m}^{k-1}\left\vert \frac{p_{m}}{P_{m}P_{m-1}}%
\tsum_{v=1}^{m}P_{v-1}\left( a_{n+v}^{m_{1}}-a_{n+v}^{m_{2}}\right)
\right\vert ^{k}\right\} ^{1/k}<\varepsilon .  \tag{2.8}
\end{equation}%
This also gives us $\left\vert a_{v}^{m_{1}}-a_{v}^{m_{2}}\right\vert
^{k}<\varepsilon /u_{0}^{1-1/k}$ for $m_{1},m_{2}>m_{0},$and all $v,$ $i.e.,$
$\left( a_{v}^{m}\right) $ is a Cauchy sequence in the set of complex
numbers $%
\mathbb{C}
.$ So it converges to a number $a_{v}$ $(v=0,1,...),$ $i.e.,$ $%
\lim_{m\rightarrow \infty }a_{v}^{m}=a_{v}.$ Now letting $m_{2}\rightarrow
\infty ,$ by $\left( 2.8\right) $ we have, for $\left\Vert
a^{m_{1}}-a\right\Vert _{\left\vert f(\overline{N}_{p}^{u})\right\vert
_{k}}<\varepsilon $ for $m_{1}>m_{0}.$ This means $\lim_{m\rightarrow \infty
}a^{m}=a.$ Further, since 
\begin{equation*}
\left\Vert a\right\Vert _{\left\vert f(\overline{N}_{p}^{u})\right\vert
_{k}}\leq \left\Vert a^{m_{1}}-a\right\Vert _{\left\vert f(\overline{N}%
_{p}^{u})\right\vert _{k}}+\left\Vert a^{m_{1}}\right\Vert _{\left\vert f(%
\overline{N}_{p}^{u})\right\vert _{k}}<\infty ,
\end{equation*}%
then, $a\in \left\vert f(\overline{N}_{p}^{u})\right\vert _{k}.$ So, $%
\left\vert f(\overline{N}_{p}^{u})\right\vert _{k}$ is a Banach space. This
completes the proof.\smallskip

We note that if $E$ is a $BK$-space such that $bs\subset E\subset \ell
_{\infty }$, then $E$ is not separable (and hence not reflexive)$\left( 
\text{see }\left[ 4\right] \right) .$ Hence the following result at once
follows from Theorem 2.2.\smallskip

\textbf{Corollary 2.5. }If $(p_{m})$ and $(u_{m})$\ are sequences of
positive numbers satisfying $\left( 2.3\right) $ and $\left( 2.4\right) $,
then $\left\vert f(\overline{N}_{p}^{u})\right\vert _{k}$ is not seperable
for $k>1.$

\section{Matrix Transformations on Space $\left\vert f(\overline{N}%
_{p})\right\vert _{k}$}

In this section we characterize certain matrix transformations on the space $%
\left\vert f(\overline{N}_{p})\right\vert _{k}.$ First we recall some
notations. Let $X,$ $Y$ be any subsets of $\omega $ and $A=(a_{nv})$ be an
infinite matrix of complex numbers. By $A(x)=(A_{n}(x)),$ we indicate the $A$%
-transform of a sequence $x=\left( x_{v}\right) $, if the series 
\begin{equation*}
A_{n}\left( x\right) =\sum\limits_{v=0}^{\infty }a_{nv}x_{v}
\end{equation*}%
are convergent for $n\geq 0.$ If $Ax\in Y,$ whenever $x\in X,$ then we say
that $A$ defines a matrix mapping from $X$ into $Y$ and denote the class of
all infinite matrices $A$ such that $A:X\rightarrow Y$ by $(X,Y).$ Also we
denote the set of all $p$-absolutely convergent series by $\ell _{k},1\leq
k<\infty ,$ $i.e.,$ 
\begin{equation*}
\ell _{k}=\left\{ x=\left( x_{v}\right) \in w:\sum\limits_{v=0}^{\infty
}\left\vert x_{v}\right\vert ^{k}<\infty \right\}
\end{equation*}%
which is a $BK$-space by respect to the norm 
\begin{equation*}
\left\Vert x\right\Vert _{\ell _{k}}=\left( \sum\limits_{v=0}^{\infty
}\left\vert x_{v}\right\vert ^{k}\right) ^{1/k}.
\end{equation*}

Also we make use of the following lemma Sar\i g\"{o}l $\left[ 14\right] $.

\textbf{Lemma 3.1. }Suppose that $A=\left( a_{nv}\right) $ is an infinite
matrix with complex numbers and $p=\left( p_{v}\right) $ is a bounded
sequence of positive numbers such that $H=\sup_{v}p_{v}~and\ $ $C=max\left\{
1,2^{H-1}\right\} .\ $Then,%
\begin{equation*}
\frac{1}{4C^{2}}U_{p}\left( A\right) \leq L_{p}(A)\leq U_{p}\left( A\right) ,
\end{equation*}%
provided that%
\begin{equation}
U_{p}\left( A\right) =\sum\limits_{v=0}^{\infty }\left(
\sum\limits_{n=0}^{\infty }\left\vert a_{nv}\right\vert \right)
^{p_{v}}<\infty
\end{equation}%
or 
\begin{equation}
L_{p}(A)=\sup \left\{ \sum\limits_{v=0}^{\infty }\left\vert
\sum\limits_{n\in N}a_{nv}\right\vert ^{p_{v}}:\text{N is a finite subset of
N}_{0}\right\} <\infty .
\end{equation}

Now we begin with first theorem given the characterization of the class $%
(\ell _{1},\left\vert f(\overline{N}_{p})\right\vert _{k}).$

\textbf{Theorem 3.2.} Let $u=(u_{n})$ be a sequence of positive numbers and
let $A=(a_{vj})$ be an infinite matrix. Then, $A\in (\ell _{1},\left\vert f(%
\overline{N}_{p})\right\vert _{k}),1\leq k<\infty ,$ if and only if%
\begin{equation}
\tsum_{m=0}^{\infty }u_{m}^{k-1}\left\vert b(m,n,j)\right\vert ^{k}<\infty 
\text{ uniformly in }n  \tag{3.1}
\end{equation}%
and 
\begin{equation}
\sup \left\{ \tsum_{m=0}^{\infty }u_{m}^{k-1}\left\vert b(m,n,j)\right\vert
^{k}:n,j\in 
\mathbb{N}
\right\} <\infty  \tag{3.2}
\end{equation}%
where 
\begin{equation*}
b(m,n,j)=\left\{ 
\begin{array}{c}
a_{nj},\ \ \ \ \ \ \ \ \ \ \ \ \ \ \ m=0 \\ 
\frac{p_{m}}{P_{m}P_{m-1}}\tsum_{v=1}^{m}P_{v-1}a_{n+v,,j},\ m\geq 1.%
\end{array}%
\right.
\end{equation*}%
\smallskip

\textbf{Proof. }Necessity. Suppose $A\in (\ell _{1},\left\vert f(\overline{N}%
_{p})\right\vert _{k}).$ Then, $A(x)\in \left\vert f(\overline{N}%
_{p})\right\vert _{k}$ for all $x\in \ell _{1},$ $i.e.,$ 
\begin{eqnarray*}
\sum_{m=0}^{\infty }u_{m}^{k-1}\left\vert F_{m,n}(A(x))\right\vert ^{k}
&=&u_{0}^{k-1}\left\vert \tsum_{j=0}^{\infty }a_{n,j}x_{j}\right\vert
^{k}+\tsum_{m=1}^{\infty }u_{m}^{k-1}\left\vert \frac{p_{m}}{P_{m}P_{m-1}}%
\tsum_{v=1}^{m}P_{v-1}\tsum_{j=0}^{\infty }a_{n+v,j}x_{j}\right\vert ^{k} \\
&=&u_{0}^{k-1}\left\vert \tsum_{j=0}^{\infty }a_{n,j}x_{j}\right\vert
^{k}+\tsum_{m=1}^{\infty }u_{m}^{k-1}\left\vert \tsum_{j=0}^{\infty }\left( 
\frac{p_{m}}{P_{m}P_{m-1}}\tsum_{v=1}^{m}P_{v-1}a_{n+v,j}\right)
x_{j}\right\vert ^{k} \\
&=&\sum_{m=0}^{\infty }u_{m}^{k-1}\left\vert \sum_{j=0}^{\infty
}b(m,n,j)x_{j}\right\vert ^{k}<\infty
\end{eqnarray*}%
uniformly in $n.$ If we put $x=e^{j}=\left( e_{v}^{j}\right) \in \ell _{1},$ 
$A(e^{j})\in \left\vert f(\overline{N}_{p})\right\vert _{k},$ where $%
e_{v}^{j}=1$ for $v=j,$ and zero otherwise, which gives that $\left(
3.2\right) $ holds. Further, \ since $\ell _{1}$ is Banach space, by the
Banach-Steinhaus theorem, $A:\ell _{1}\rightarrow 
\mathbb{C}
$ is a continuous linear map. So, for fixed $n$ and $s,$ 
\begin{equation*}
q_{sn}(x)=\left( \sum_{m=0}^{s}u_{m}^{k-1}\left\vert
F_{m,n}(A(x))\right\vert ^{k}\right) ^{1/k}
\end{equation*}%
is a continuous seminorm on $\ell _{1},$ which implies that $%
\lim_{s\rightarrow \infty }q_{sn}(x)=q_{n}(x)$ is a continuous seminorm, or,
equivalently, there exists a constant $K$ such that 
\begin{equation}
q_{n}(x)=\left( \sum_{m=0}^{\infty }u_{m}^{k-1}\left\vert
F_{m,n}(A(x))\right\vert ^{k}\right) ^{1/k}\leq K\left\Vert x\right\Vert
_{\ell _{1}}  \tag{3.3}
\end{equation}%
for every $x\in \ell _{1}.$ Applying $\left( 3.3\right) $ with $%
x=e^{j}=\left( e_{v}^{j}\right) \in \ell _{1}$\ we have, for all $j,n\geq 0,$%
\begin{equation*}
\left( \sum_{m=0}^{\infty }u_{m}^{k-1}\left\vert b(m,n,j\right\vert
^{k}\right) ^{1/k}\leq K,
\end{equation*}%
which gives $\left( 3.2\right) .$

Sufficiency. Suppose $\left( 3.1\right) $ and $\left( 3.2\right) $ hold.
Given $x\in \ell _{1}.$ Then, we should show $A(x)\in \left\vert f(\overline{%
N}_{p})\right\vert _{k}.$ For this, it is enough to prove that 
\begin{equation*}
\tsum_{m=l}^{\infty }u_{m}^{k-1}\left\vert F_{m,n}(A(x)))\right\vert
^{k}\rightarrow 0\text{ as }l\rightarrow \infty ,\text{uniformly in }n.
\end{equation*}%
By applying Minkowski's inequality we get 
\begin{eqnarray}
\left( \sum_{m=l}^{\infty }u_{m}^{k-1}\left\vert F_{m,n}(A(x))\right\vert
^{k}\right) ^{1/k} &=&\left( \tsum_{m=M}^{\infty }u_{m}^{k-1}\left\vert
\tsum_{j=0}^{\infty }b(m,n,j)x_{j}\right\vert ^{k}\right) ^{1/k}  \notag \\
&\leq &\tsum_{j=0}^{\infty }\left\vert x_{j}\right\vert \left(
\tsum_{m=l}^{\infty }u_{m}^{k-1}\left\vert b(m,n,j)\right\vert ^{k}\right)
^{1/k}  \TCItag{3.4} \\
&=&\tsum_{j=0}^{\infty }\left\vert x_{j}\right\vert R(l,n,j)  \notag
\end{eqnarray}%
where 
\begin{equation*}
R(l,n,j)=\left( \tsum_{m=l}^{\infty }u_{m}^{k-1}\left\vert
b(m,n,j)\right\vert ^{k}\right) ^{1/k}.
\end{equation*}%
On the other hand, by $\left( 3.2\right) ,$ since for all $l,n,j\geq 0,$ 
\begin{equation*}
R(l,n,j)\left\vert x_{j}\right\vert \leq \left\vert x_{j}\right\vert \sup
\left\{ R(0,n,j):n,j\in 
\mathbb{N}
\right\}
\end{equation*}%
it follows that 
\begin{equation*}
\tsum_{j=0}^{\infty }\left\vert x_{j}\right\vert R(l,n,j)<\infty \text{
uniformly in }n,l
\end{equation*}%
Hence, for every $\varepsilon >0,$ there exists an integer $j_{0}$ such
that, for all $l$ and $n,$ 
\begin{equation*}
\tsum_{j=j_{0}}^{\infty }\left\vert x_{j}\right\vert R(l,n,j)<\frac{%
\varepsilon }{2}.
\end{equation*}%
Also, by $\left( 3.1\right) ,$ since%
\begin{equation*}
R(l,n,j)\rightarrow 0\text{ as }l\rightarrow \infty \text{ uniformly in }n,
\end{equation*}%
there exists integer $l_{0}$ so that, for $l\geq l_{0}$ and all $n,$%
\begin{equation*}
\tsum_{j=0}^{j_{0}-1}\left\vert x_{j}\right\vert R(l,n,j)<\frac{\varepsilon 
}{2}.
\end{equation*}%
So we have%
\begin{equation*}
\tsum_{j=0}^{\infty }\left\vert x_{j}\right\vert R(l,n,j)<\varepsilon
\end{equation*}%
which implies, by $\left( 3.4\right) ,$ 
\begin{equation*}
\sum_{m=l}^{\infty }u_{m}^{k-1}\left\vert F_{m,n}(A(x))\right\vert
^{k}\rightarrow 0\text{ as }l\rightarrow \infty \text{ informly in }n.
\end{equation*}%
This completes the proof. \smallskip

In the special case $p_{m}=u_{m}=1$ for all $m\geq 0,$ $\left\vert f(%
\overline{N}_{p})\right\vert =\widehat{\ell }_{k}$ and so the following
result follows from Theorem 3.1.

\textbf{Corollary 3.3. }$A\in (\ell _{1},\widehat{\ell }_{k}),1\leq k<\infty
,$ if and only if%
\begin{equation*}
\tsum_{m=0}^{\infty }\left\vert b(m,n,j)\right\vert ^{k}<\infty \text{
uniformly in }n
\end{equation*}%
and 
\begin{equation*}
\sup \left\{ \tsum_{m=0}^{\infty }\left\vert b(m,n,j)\right\vert ^{k}:n,j\in 
\mathbb{N}
\right\} <\infty
\end{equation*}%
where 
\begin{equation}
b(m,n,j)=\left\{ 
\begin{array}{c}
a_{nj},\ \ \ \ \ \ \ \ \ \ \ \ \ \ \ m=0 \\ 
\frac{1}{m(m+1)}\tsum_{v=1}^{m}va_{n+v,,j},\ m\geq 1.%
\end{array}%
\right.  \tag{3.5}
\end{equation}

\textbf{Theorem 3.4.} Let $u=(u_{n})$ be a sequence of positive numbers and
let $A=(a_{vj})$ be an infinite matrix . Then, $A\in (c,\left\vert f(%
\overline{N}_{p})\right\vert _{k}),1\leq k<\infty ,$ if and only if%
\begin{equation}
B=\sup \left\{ \tsum_{m=0}^{\infty }u_{m}^{k-1}\left( \tsum_{j=0}^{\infty
}\left\vert b(m,n,j)\right\vert \right) ^{k}:n\in 
\mathbb{N}
\right\} <\infty ,\text{ }  \tag{3.6}
\end{equation}%
\begin{equation}
\tsum_{m=0}^{\infty }u_{m}^{k-1}\left\vert \tsum_{j=0}^{\infty
}b(m,n,j)\right\vert ^{k}<\infty \text{ uniformly in }n,  \tag{3.7}
\end{equation}%
\begin{equation}
\tsum_{m=0}^{\infty }u_{m}^{k-1}\left\vert b(m,n,j)\right\vert ^{k}<\infty 
\text{ uniformly in }n.  \tag{3.8}
\end{equation}

\textbf{Proof.} Necessity. Let $A\in (c,\left\vert f(\overline{N}%
_{p})\right\vert _{k}).$ Then, $A(x)\in \left\vert f(\overline{N}%
_{p})\right\vert _{k}$ for every $x\in c.$ Now, let $e,e^{\left( j\right)
}\in c.$ Then, we have $\left( 3.7\right) $ and $\left( 3.8\right) ,$
respectively, where $e=(1,1,...).$ Also, it follows as in the proof of
Theorem 3.2 that 
\begin{equation}
\left( \sum_{m=0}^{\infty }u_{m}^{k-1}\left\vert F_{m,n}(A(x))\right\vert
^{k}\right) ^{1/k}\leq K\left\Vert x\right\Vert _{\infty }.  \tag{3.9}
\end{equation}%
Let $N$ be arbitrary finite set of natural numbers and define a sequence $x$
by 
\begin{equation}
x_{j}=\left\{ 
\begin{array}{c}
1,\ j\in N \\ 
0,\text{\ }j\notin N,%
\end{array}%
\right.  \tag{3.10}
\end{equation}%
then $x\in c$ and $\left\Vert x\right\Vert =1.$ Applying $\left( 3.9\right) $
with this sequence $\left( 3.10\right) ,$ we have 
\begin{equation}
\left\{ \tsum_{m=0}^{\infty }u_{m}^{k-1}\left( \left\vert \tsum_{j\in
N}b(m,n,j)\right\vert \right) ^{k}\right\} ^{1/k}\leq K.  \tag{3.11}
\end{equation}%
Hence, it is seen from Lemma 3.1 together with $p_{v}=1$ for all $v$ that $%
\left( 3.11\right) $ is equivalent to $\left( 3.6\right) .$

Sufficiency. Suppose that $\left( 3.6\right) ,\left( 3.7\right) $ and $%
\left( 3.8\right) $ hold. Given $x\in c$ and say $\lim_{j}x_{j}=\beta .$
Then, by $\left( 3.7\right) ,$ as in Theorem 3.2, 
\begin{equation*}
\tsum_{m=0}^{\infty }u_{m}^{k-1}\left\vert F_{mn}(A(x))\right\vert
^{k}=\tsum_{m=0}^{\infty }u_{m}^{k-1}\left\vert \tsum_{j=0}^{\infty
}b(m,n,j)x_{j}\right\vert ^{k}<\infty .
\end{equation*}%
Now it is enough to show that the tail of this series tends to zero
uniformly in $n.$ To see that we write 
\begin{equation*}
\tsum_{m=M}^{\infty }u_{m}^{k-1}\left\vert F_{mn}(A(x))\right\vert
^{k}=\tsum_{m=M}^{\infty }u_{m}^{k-1}\left\vert \tsum_{j=0}^{\infty }\beta
b(m,n,j)+\tsum_{j=0}^{\infty }b(m,n,j)\left( x_{j}-\beta \right) \right\vert
^{k}
\end{equation*}%
\begin{eqnarray*}
&=&O(1)\left\{ \tsum_{m=M}^{\infty }u_{m}^{k-1}\left\vert
\tsum_{j=0}^{\infty }\beta b(m,n,j)\right\vert ^{k}+\tsum_{m=M}^{\infty
}u_{m}^{k-1}\left\vert \tsum_{j=0}^{\infty }b(m,n,j)\left( x_{j}-\beta
\right) \right\vert ^{k}\right\} \\
&=&O(1)\left\{ F_{M,n}^{1}+F_{M,n}^{2}(x^{\prime })\right\} ,\text{ }say.
\end{eqnarray*}%
It is clear from $\left( 3.7\right) $ that $F_{M,n}^{1}\rightarrow 0$ as $%
M\rightarrow \infty $ uniformly in $n$. On the other hand, since $%
x_{j}\rightarrow \beta ,$ for any $\varepsilon >0,$ there exists an integer $%
j_{0}$ such that 
\begin{equation*}
\left\vert x_{j}-\beta \right\vert <\frac{\varepsilon ^{1/k}}{2B}\text{ \
for }j\geq j_{0}
\end{equation*}%
which gives us, by $\left( 3.6\right) ,$ for all $n\geq 0,$ 
\begin{eqnarray*}
F_{M,n}^{2}(x^{\prime }) &=&\tsum_{m=M}^{\infty }u_{m}^{k-1}\left\vert
\tsum_{j=0}^{\infty }b(m,n,j)\left( x_{j}-\beta \right) \right\vert ^{k} \\
&=&O(1)\tsum_{m=M}^{\infty }u_{m}^{k-1}\left\{ \left\vert
\tsum_{j=0}^{j_{0}-1}b(m,n,j)\left( x_{j}-\beta \right) \right\vert
^{k}+\left\vert \tsum_{j=j_{0}}^{\infty }b(m,n,j)\left( x_{j}-\beta \right)
\right\vert ^{k}\right\} \\
&=&O(1)\tsum_{m=M}^{\infty }u_{m}^{k-1}\left(
\tsum_{j=0}^{j_{0}-1}\left\vert b(m,n,j)\right\vert \right) ^{k}+\frac{%
\varepsilon }{2}.
\end{eqnarray*}%
By $\left( 3.8\right) ,$ the first term of the equality is smaller that $%
\varepsilon /2$ for suffiently large $M$ and all $n.$ This means $%
F_{M,n}^{2}\rightarrow 0$ as $M\rightarrow \infty $ uniformly in $n$. Hence,
the theorem is established.

For $p_{m}=u_{m}=1,$ Theorem 3.4 is reduced to the following result.

\textbf{Corollary 3.5.} Let $A=(a_{vj})$ be an infinite matrix and $\left(
b(m,n,j)\right) $ be as in $\left( 3.5\right) .$ Then, $A\in (c,\widehat{%
\ell }_{k}),1\leq k<\infty ,$ if and only if%
\begin{equation*}
B=\sup \left\{ \tsum_{m=0}^{\infty }\left( \tsum_{j=0}^{\infty }\left\vert
b(m,n,j)\right\vert \right) ^{k}:n\in 
\mathbb{N}
\right\} <\infty \text{ }
\end{equation*}%
\begin{equation*}
\tsum_{m=0}^{\infty }\left\vert \tsum_{j=0}^{\infty }b(m,n,j)\right\vert
^{k}<\infty \text{ uniformly in }n
\end{equation*}%
\begin{equation*}
\tsum_{m=0}^{\infty }\left\vert b(m,n,j)\right\vert ^{k}<\infty \text{
uniformly in }n
\end{equation*}

\end{document}